\documentclass {article}
\usepackage{authblk}
\usepackage{etex}
\usepackage[utf8]{inputenc}
\usepackage[english]{babel}
\usepackage{amsmath}
\usepackage{amsthm}
\usepackage{amssymb}
\usepackage{titlesec}
\usepackage{color}
\usepackage[color,matrix,arrow]{xy}
\usepackage{amsgen}
\usepackage{amstext}
\usepackage{amsbsy}
\usepackage{amsopn}
\usepackage{amsfonts}
\usepackage{eepic}
\usepackage{graphicx}
\usepackage{epsf}
\usepackage{pstricks}
\usepackage{float}
\usepackage{enumerate}
\usepackage{eqnarray}
\usepackage{faktor}
\usepackage{ mathrsfs }
\xyoption{all}

\usepackage{pgf}
\usepackage{tikz-cd}
\usetikzlibrary{automata, arrows.meta, positioning,decorations.pathmorphing}

\newcommand{\footrecall}[1]{
} 
\usetikzlibrary{snakes,shapes,arrows,automata}

\titleformat*{\section}{\large\bfseries}
\titleformat*{\subsection}{\normalsize \bfseries}

\newcommand{\Z}{\mathbb{Z}}

\newcommand{\End}{\text{End}}
\newcommand{\Aut}{\text{Aut}}
\newcommand{\Rat}{\text{Rat}}
\newcommand{\Rec}{\text{Rec}}
\newcommand{\Alg}{\text{Alg}}
\newcommand{\CF}{\text{CF}}

\newcommand{\mc}{\mathcal}

\theoremstyle{definition}
\newtheorem{theorem}{Theorem}[section]
\newtheorem{corollary}[theorem]{Corollary}

\newtheorem{proposition}[theorem]{Proposition}

\newtheorem{problem}[theorem]{Problem}
\newtheorem{lemma}[theorem]{Lemma}

\begin{document}

\title{Computing centralizers in [f.g. free]-by-cyclic groups}
\author{Andr\'e Carvalho
\\

Center for Mathematics and Applications (NOVA Math), NOVA SST

2829–516 Caparica, Portugal

andrecruzcarvalho@gmail.com
}
\maketitle

\begin{abstract}
We prove that centralizers of elements in [f.g. free]-by-cyclic groups are computable. As a corollary we get that, given two conjugate elements in a [f.g. free]-by-cyclic group, the set of conjugators can be computed and that the conjugacy problem with context-free constraints is decidable. In the end, we pose several problems arising naturally from this work.
\end{abstract}
\section{Introduction}
Two elements $x,y\in G$ are said to be \emph{conjugate} if there is some $z\in G$ such that $x=z^{-1}yz$. Conjugacy is an equivalence relation and we write $x\sim y$ if $x$ and $y$ are conjugate. The \emph{conjugacy problem}  is an important problem in combinatorial group theory and consists on deciding whether two elements in a group $G$  are {conjugate} or not. It is undecidable in general \cite{[Mil92]}  and has been studied by many in different classes of groups (see, for example \cite{[Gro87],[OS06],[BMMV06],[LS11],[Log22]}).

Some variations and questions adjacent to the conjugacy problem have also been studied. In \cite{[LS11]}, Ladra and Silva study the \emph{generalized conjugacy problem}, which is the problem of deciding, given an element $g \in G$ and a subset $K\subseteq G$, whether $g$ has a conjugate in $K$ and prove that this is decidable in virtually free groups when the subset $K$ is rational. More than that, they prove that the set of solutions of the equation $x^{-1}gx\in K$ is rational and effectively constructible in a virtually free group. This allows the solution of the \emph{constrained generalized conjugacy problem} when the constraints have good intersection properties with rational subsets. This problem consists on deciding whether $g$ has a conjugate in $K$ with the conjugator belonging to some constraint subset.  Since it is decidable whether a context-free and a rational subset intersect, it follows that the generalized conjugacy problem for rational subsets of virtually free groups with context-free subsets as constraints is decidable. We remark that the class of context-free subsets of virtually free groups contains the class of rational subsets (see \cite{[Her91]}) and so, in particular, we can decide this problem when the constraint is a finitely generated subgroup. 
Even though this is not relevant in the virtually free group case since finitely generated subgroups of virtually free groups are again virtually free, we remark that solving the conjugacy problem with finitely generated subgroups as constraints in a given group $G$ yields a solution of the conjugacy problem in any finitely generated subgroup $H$ of $G$, since two elements $h_1,h_2\in H$ are conjugate in $H$ if and only if they are conjugate in $G$ with a conjugator belonging to $H$.

Given an endomorphism $\phi\in \End(G)$, two elements $x,y\in G$ are $\phi$-twisted conjugate if there is some $z\in G$ such that $x=(z^{-1}\phi) y z$, in which case we write $x\sim_{\phi}y$. This is again an equivalence relation. Brinkmann's conjugacy problem was introduced in \cite{[Bri10]}, where it was solved in the affirmative for automorphisms of the free group, and consists on deciding, given an automorphism $\phi$ and elements $x,y\in G$, whether there is some integer $k$ such that $x\phi^k\sim y$.

In \cite{[BMMV06]}, Bogopolski, Martino, Maslakova and Ventura proved that the conjugacy problem is decidable in [f.g. free]-by-cyclic groups by proving that it can be reduced to the twisted conjugacy problem and Brinkmann's conjugacy problem for automorphisms of the free group. They solve the twisted conjugacy problem and use Brinkmann's result from \cite{[Bri10]} to finish the proof. This approached was followed in \cite{[Car23b]} for the generalized versions of these problems. 
In \cite{[Log22]}, Logan generalized the ideas in \cite{[BMMV06]} and proved that the conjugacy problem is decidable in ascending HNN-extensions of the free group by reducing it to non-bijective versions of Brinkmann's conjugacy problem and the twisted conjugacy problem and proving decidability of these algorithmic questions. In \cite{[Log22]}, Logan highlights the fact that \emph{many} (in a precise sense) one-relator groups appear as subgroups of ascending HNN-extensions of a free group, so solving the conjugacy problem for ascending HNN-extensions with finitely generated subgroups as constraints would yield a solution to the conjugacy problem in many one-relator groups, which is an important open problem.

The main goal of this paper is to prove that centralizers of elements of [f.g. free]-by-cyclic groups are computable. If the cyclic group is finite, then our group is virtually free, so this is well-known. For this reason, we will focus on the [f.g. free]-by-$\Z$ case. 
\newtheorem*{thm centralizers}{Theorem \ref{thm centralizers}}
\begin{thm centralizers}
There is an algorithm taking as input  an automorphism $\phi\in \Aut(F_n)$ and an element $t^ax\in F_n\rtimes_\phi \Z$ that outputs a finite set of generators for $C_{F_n\rtimes_\phi \Z}(t^ax)$.
\end{thm centralizers}

Naturally, if two elements $x,y\in G$ are conjugate, the set of conjugators is a coset of $C_G(x)$ and, since the conjugacy problem is decidable in [f.g. free]-by-$\Z$ groups, then we have the following corollary, in the same spirit as the result in \cite{[LS11]}.

\newtheorem*{cor: solseq}{Corollary \ref{cor: solseq}}
\begin{cor: solseq}
Let $\phi\in \Aut(F_n)$ and $t^au, \, t^bv\in F_n\rtimes_\phi\Z$. The set of solutions  to the equation $x^{-1}(t^au)x=t^bv$ is rational and effectively constructible.
\end{cor: solseq}

This corollary clearly implies that the solution set of the equation $x^{-1}(t^au)x\in K$ for a finite $K$ is rational and computable since it is a finite union of rational and computable subsets. Example 5.3 in \cite{[LS11]} shows that we cannot replace finite by rational, as in this case the solution set is not rational. 

Since   context-free languages are closed under the intersection with rational languages and we can decide its emptiness, we also have the following corollary:

\newtheorem*{cor: cf constraint}{Corollary \ref{cor: cf constraint}}
\begin{cor: cf constraint}
The conjugacy problem with context-free constraints is decidable for [f.g. free]-by-cyclic groups, i.e., there is an algorithm taking as input an two elements $t^ax,\, t^by \in F_n\rtimes_{\phi}\Z$ and a context-free grammar generating a language $L$ such that $L=K\pi^{-1}$ for some $K\subseteq F_n\rtimes_\phi \Z$, that decides whether there is an element $t^cz\in K$ such that $t^ax=(t^cz)^{-1}(t^by) t^c z$. 
\end{cor: cf constraint}

\section{Preliminaries}
We will now present the basic definitions and results on rational, algebraic and context-free subsets of groups. For more detail, the reader is referred to \cite{[Ber79]} and \cite{[BS21]}.

Let $G=\langle A\rangle$ be a finitely generated group, $A$ be a finite generating set, $\tilde A=A\cup A^{-1}$ and $\pi:\tilde A^*\to G$ be the canonical (surjective) homomorphism. 

A subset $K\subseteq G$ is said to be \emph{rational} if there is some rational language $L\subseteq \tilde A^*$ such that $L\pi=K$ and \emph{recognizable} if $K\pi^{-1}$ is rational. 
We will denote by $\Rat(G)$ and $\Rec(G)$ the class of rational and recognizable subsets of $G$, respectively. Rational subsets generalize the notion of finitely generated subgroups.

\begin{theorem}[\cite{[Ber79]}, Theorem III.2.7]
\label{AnisimovSeifert}
Let $H$ be a subgroup of a group $G$. Then $H\in \Rat(G)$ if and only if $H$ is finitely generated.
\end{theorem}

We remark that, as long as we can test membership in the subset $H$, the above theorem is constructive, in the sense that if we have an automaton recognizing a language $L$ such that $L\pi$ is a subgroup, we can compute a set of generators for the subgroup $H$. That can be seen by following the proof in \cite[Theorem 3.1]{[BS21]}.

A natural generalization of these concepts relates to the class of context-free languages. 
A subset $K\subseteq G$ is said to be \emph{algebraic} if there is some context-free language $L\subseteq \tilde A^*$ such that $L\pi=K$ and \emph{context-free} if $K\pi^{-1}$ is context-free. 
We will denote by $\Alg(G)$ and $\CF(G)$ the class of algebraic and context-free subsets of $G$, respectively.  It follows from \cite[Lemma 2.1]{[Her91]} that $\CF(G)$ and $\Alg(G)$ do not depend on the alphabet $A$ or the surjective homomorphism $\pi$. We follow the terminology in \cite{[Her91], [Her92], [Car23c]}. However, for example in \cite{[CEL23]}, algebraic subsets are called context-free subsets and context-free subsets are called recognisably context-free.

It is obvious from the definitions that $\Rec(G)$, $\Rat(G)$, $\CF(G)$ and $\Alg(G)$ are closed under union, since both rational and context-free languages are closed under union. The intersection case is distinct: from the fact that rational languages are closed under intersection, it follows that $\Rec(G)$ must be closed under intersection too. However $\Rat(G)$, $\Alg(G)$ and $\CF(G)$ might not be. We will also use the fact that the class of rational subsets is closed under inversion.

For a finitely generated group $G$, it is immediate from the definitions that $$\Rec(G)\subseteq \CF(G) \subseteq \Alg(G)$$ and that $$\Rec(G)\subseteq \Rat(G) \subseteq \Alg(G).$$ It is proved in \cite{[Her91]} that 
$$\CF(G)=\Alg(G) \Leftrightarrow \CF(G)=\Rat(G) \Leftrightarrow \text{ G is virtually cyclic.}$$

However, there is no general inclusion between $\Rat(G)$ and $\CF(G)$. For example, if $G$ is virtually abelian, then $\CF(G)\subseteq \Alg(G)= \Rat(G)$ (and the inclusion is strict if the group is not virtually cyclic) and if the group is virtually free, then $\Rat(G)\subseteq \CF(G)$ (see \cite[Lemma 4.2]{[Her91]}). 
In the case of the free group $F_n$ of rank $n\geq 1$, Herbst proves in  \cite[Lemma 4.6]{[Her91]} an analogue of Benois' Theorem for context-free subsets, proving that for a subset $K\subseteq F_n$, then $K\in \CF(F_n)$ if and only if the set of reduced words representing elements of $K$ is context-free.

\section{The main result}

Let $G$ be a group, $\phi\in \Aut(G)$, $x\in G$ and $a\in \Z$. We define $\mc E_{a,x,\phi}=\{k\in \Z\mid x\phi^k\sim_{\phi^a} x\}$. Usually we omit the subscripts $\phi$ and $x$ as these will be clear from the context.
\begin{lemma}\label{lem: cyclic}
Let $G$ be a group, $\phi\in \Aut(G)$, $x\in G$ and $a\in \Z$. Then, either $\mc E_{a,x,\phi}=\{0\}$ or there  is some $b$ dividing $a$ such that $$\mc E_{a,x,\phi}=b\Z.$$
\end{lemma}
\noindent\textit{Proof.}  We have that $0\in \mc E_{a}$, because $x\phi^0=x=(1^{-1}\phi^a)x1$. Suppose that $\mc E_{a}\neq \{0\}$. Since $x\phi^a=(x\phi^a)xx^{-1}$, then $x\sim_{\phi^a}x\phi^a$, and so $a\in \mc E_{a}$. Now, let $b=\min\{|k|\mid k\in \mc E_{a}\setminus\{0\}\}$. We will prove that $\mc E_{a}=b\Z$, which in particular implies that $b\mid a$, as $a\in \mc E_{a}$.
We start by showing that $k\in \mc E_{a}$ if and only if $-k\in \mc E_{a}$. This follows from the fact that
$$x\phi^b=(y^{-1}\phi^a)xy\iff x=(y^{-1}\phi^{a-b})(x\phi^{-b})(y\phi^{-b})\iff x\phi^{-b}=(y\phi^{-b})\phi^ax(y^{-1}\phi^{-b}).$$
Now we show that $\mc E_{a}$ is closed under addition, from where it follows that $\mc E_{a}$ is a subgroup of $\Z$, and so, cyclic.
Let $k_1,k_2\in\mc E_{a}$. Then, there are $y,z\in G$ such that  $x\phi^{k_1}=(y^{-1}\phi^a)xy$ and $x\phi^{k_2}=(z^{-1}\phi^a)xz$. Hence,
\begin{align*}
x\phi^{k_1+k_2}=&x\phi^{k_1}\phi^{k_2}\\
=&((y^{-1}\phi^a)xy)\phi^{k_2}\\
=&(y^{-1}\phi^{k_2})\phi^a(x\phi^{k_2})(y\phi^{k_2})\\
=&(y^{-1}\phi^{k_2})\phi^a((z^{-1}\phi^a)xz)(y\phi^{k_2})\\
=&((y^{-1}\phi^{k_2})z^{-1})\phi^ax(z(y\phi^{k_2})).
\end{align*}
\qed\\
Notice that, since $a$ necessarily belongs to $\mc E_{a,x,\phi}$, then $\mc E_a=0$ implies $a=0$.
We will denote the minimal element $b$ from the lemma above by $e_a$.

\begin{lemma}\label{centralizerfree}
Let $G$ be a group, $\phi\in \Aut(G)$ and $x\in G$. We have that:
\begin{enumerate}[i.]
\item if there is  $k\neq 0$  such that $x\phi^k\sim x$, then  
$$\faktor{C_{G\rtimes_\phi \Z}(x)}{C_G(x)}=\langle (t^{e_0}z)C_G(x)\rangle,$$
 where $z\in G$ is such that $x=z^{-1}(x\phi^{e_0})z$. In particular $$C_{G\rtimes_\phi \Z}(x)=\langle t^{e_0}z, C_{G}(x)\rangle;$$
\item if there is  no  $k\neq 0$ such that $x\phi^k \sim x$, then 
$$C_{G\rtimes_\phi \Z}(x)= C_{G}(x).$$
\end{enumerate}
\end{lemma}
\noindent\textit{Proof.}  We have that, for $t^b y\in G\rtimes_\phi\Z$,
\begin{align}\label{c1}
x (t^b y)=(t^b y)x\iff(y^{-1}t^{-b})x (t^b y)= x \iff y^{-1}(x\phi^b)y=x.
\end{align}
If there is no $k\neq 0$ such that $x\phi^k \sim x$, then, by (\ref{c1}), $t^by\in C_{G\rtimes_\phi \Z}$ implies that $b=0$ and $y\in C_G(x)$. So, in this case, we have that
$C_{G\rtimes_\phi \Z}(x)=C_{G}(x)$.

Assume now that there is such a $k$. Clearly, $C_G(x)\trianglelefteq C_{G\rtimes_\phi \Z}(x)$ and $t^{e_0}z\in C_{G\rtimes_\phi \Z}(x)$ by (\ref{c1}). We have that, if  $t^ay\in C_{G\rtimes_\phi \Z}(x)$, then $a\in e_0\Z$, by Lemma \ref{lem: cyclic}, so $a=\lambda e_0$, for some $\lambda\in \Z$. This means that $(t^{e_0}z)^{-\lambda}(t^ay)\in G\cap  C_{G\rtimes_\phi \Z}(x)= C_{G}(x)$, and so $(t^ay)\in (t^{e_0}z)^{\lambda}C_G(x)= ((t^{e_0}z)C_G(x))^{\lambda}$.
\qed\\

 Now, we define $$\mc C_{k,a,x,\phi}:=\{y\in G\mid x=(y^{-1}\phi^a)(x\phi^k)y\}.$$ Notice that $\mc C_{k,a,x,\phi}\neq \emptyset$ if and only if $k\in \mc E_{a,x,\phi}$.
Again, we will typically write $\mc C_{k}$, since $a$, $x$ and $\phi$ will be clear from context.
\begin{proposition}\label{centralizernonfree}
Let $G$ be a group, $\phi\in \Aut(G)$, $a\in \Z$, $x\in G$ and $k\in \Z$. If $\mc E_{a,x,\phi}\neq \{0\}$, then
$$\mc C_{(k+1)e_a,a,x,\phi}= (\mc C_{ke_a,a,x,\phi})\phi^{e_a}\mc C_{e_a,a,x,\phi}= (\mc C_{ke_a,a,x,\phi})\phi^{e_a}z,$$
for all $z\in  \mc C_{e_a,a,x,\phi}$.
\end{proposition}
\noindent\textit{Proof.}  Let $y\in\mc C_{(k+1)e_a}$. Then 
\begin{align}\label{yxk1ea}
x=(y^{-1}\phi^a)(x\phi^{(k+1)e_a})y.
\end{align}
Let $z\in \mc C_{e_a}$, i.e., $z\in G$ such that 
\begin{align}\label{zxea}
x=(z^{-1}\phi^a)(x\phi^{e_a})z,
\end{align}
and so that
\begin{align}\label{x-eaz}
x\phi^{-e_a}=(z^{-1}\phi^{a-e_a})x(z\phi^{-e_a}).
\end{align}

From (\ref{yxk1ea}), we deduce that 
$x\phi^{-e_a}=(y^{-1}\phi^{a-e_a})(x\phi^{ke_a})(y\phi^{-e_a})$, which, combined with (\ref{x-eaz}), yields
$$(z^{-1}\phi^{a-e_a})x(z\phi^{-e_a})=(y^{-1}\phi^{a-e_a})(x\phi^{ke_a})(y\phi^{-e_a}),$$
or, equivalently, 
$$x=((zy^{-1})\phi^{-e_a})\phi^a(x\phi^{ke_a})(yz^{-1})\phi^{-e_a},$$
i.e., $(yz^{-1})\phi^{-e_a}\in \mc C_{ke_a}$. Hence, $y\in   (\mc C_{ke_a})\phi^{e_a}z$, and it follows that $$\mc C_{(k+1)e_a}\subseteq  (\mc C_{ke_a})\phi^{e_a}z\subseteq   (\mc C_{ke_a})\phi^{e_a}\mc  C_{e_a}.$$

Now, let $y\in \mc C_{ke_a}$ and $z\in \mc C_{e_a}$. Then $x=(y^{-1}\phi^a)(x\phi^{ke_a})y$ and $x=(z^{-1}\phi^a)(x\phi^{e_a})z$.
From the first condition, we have that $x\phi^{e_a}=(y^{-1}\phi^{a+e_a})(x\phi^{(k+1)e_a})(y\phi^{e_a})$, which, from the second condition yields that 
$$(z\phi^a)xz^{-1}=(y^{-1}\phi^{a+e_a})(x\phi^{(k+1)e_a})(y\phi^{e_a}),$$
or, equivalently,
$$x=(z^{-1}(y^{-1}\phi^{e_a}))\phi^a(x\phi^{(k+1)e_a})((y\phi^{e_a})z),$$
i.e., $(y\phi^{e_a})z\in \mc C_{(k+1)e_a}$
Hence, $  (\mc C_{ke_a})\phi^{e_a}\mc C_{e_a}\subseteq \mc C_{(k+1)e_a}$. Therefore, we have that 
$$\mc C_{(k+1)e_a,a,x,\phi}= (\mc C_{ke_a,a,x,\phi})\phi^{e_a}\mc C_{e_a,a,x,\phi}= (\mc C_{ke_a,a,x,\phi})\phi^{e_a}z,$$
for all $z\in  \mc C_{e_a,a,x,\phi}$.
\qed\\

\begin{theorem}\label{thm centralizers}
There is an algorithm taking as input  an automorphism $\phi\in \Aut(F_n)$ and an element $t^ax\in F_n\rtimes_\phi \Z$ that outputs a finite set of generators for $C_{F_n\rtimes_\phi \Z}(t^ax)$.
\end{theorem}
\noindent\textit{Proof.}
We start with the case where $a=0$. We start by computing $u\in F_n$ such that  $C_{F_n}(x)=\langle u\rangle$. Then we decide if there is some $k\neq 0$ such that $x\phi^k \sim x$, using \cite[Lemma 3.2]{[Log22]} with input $(\phi,x\phi,x)$. If there is, we compute the minimal $k>0$ such that $x\phi^k\sim x$ by solving the conjugacy problem until we find a positive answer and compute a conjugator $z$ in $F_n$. Lemma \ref{centralizerfree} yields that $C_{F_n\rtimes_\phi \Z}(x)=\langle t^kz, u\rangle$. If there is no such $k$, then Lemma \ref{centralizerfree} yields that $C_{F_n\rtimes_\phi \Z}(x)=\langle u\rangle$.

Now we have the case where $a\neq 0$. We have that 

$$(y^{-1}t^{-b})(t^ax) (t^b y)=t^a x \iff t^a (y^{-1}\phi^a)(x\phi^b)y=t^ax\iff (y^{-1}\phi^a)(x\phi^b)y=x.$$
So $$t^by\in C_{F_n\rtimes_\phi \Z}(t^ax)\iff y\in \mc C_{b,a,x,\phi}$$
and, in particular, $b\in \mc E_{a,x,\phi}$.
From Lemma \ref{lem: cyclic} follows that  $\mc E_{a,x,\phi}=e_a\Z$ (notice that $a\in \mc E_{a,x,\phi}$ and $a\neq 0$). Also, since the twisted conjugacy problem is decidable for automorphisms of free groups by \cite[Theorem 1.5]{[BMMV06]}, then $e_a$ is computable.

Hence, $$C_{F_n\rtimes_\phi \Z}(t^ax)=\bigcup_{k\in \Z} t^{ke_a}\mc C_{ke_a,a,x,\phi}$$

By Proposition \ref{centralizernonfree}, for all $k\in \Z$, $$\mc C_{(k+1)e_a,a,x,\phi}= \mc C_{ke_a,a,x,\phi}\phi^{e_a}\mc C_{e_a,a,x,\phi}.$$  It follows by induction that, for $k\geq 0$,  $$\mc C_{ke_a,a,x,\phi}= (\mc C_{e_a,a,x,\phi}\phi^{ke_a})(\mc C_{e_a,a,x,\phi}\phi^{(k-1)e_a})\cdots (\mc C_{e_a,a,x,\phi}\phi^{e_a})\mc C_{e_a,a,x,\phi}.$$

By  \cite[Proposition 5.7]{[LS11]}, $\mc C_{e_a,a,x,\phi}$ is rational and effectively constructible.
It follows that that
\begin{align*}
\bigcup_{k>0} t^{ke_a}\mc C_{ke_a,a,x,\phi}=\, & \bigcup_{k>0} t^{ke_a} \mc C_{e_a,a,x,\phi}\phi^{ke_a}\mc C_{e_a,a,x,\phi}\phi^{(k-1)e_a}\cdots \mc C_{e_a,a,x,\phi}\\
=\,& \bigcup_{k>0}(\mc C_{e_a,a,x,\phi}t^{e_a})^k\mc C_{e_a,a,x,\phi}\\
=\,& (\mc C_{e_a,a,x,\phi}t^{e_a})^+\mc C_{e_a,a,x,\phi}
\end{align*}
 and so   $\mc S=\bigcup_{k>0} t^{ke_a}\mc C_{ke_a,a,x,\phi}$ is also rational and effectively constructible.
 Again, the set of elements with negative integral part is the set of inverses of elements in $\mc S$. Since rational subsets are effectively closed under taking inverses and, again by  \cite[Proposition 5.7]{[LS11]},   $\mc C_{0,a,x,\phi}=\{y\in G\mid (y^{-1}\phi^a) x y=x\}$ is rational and computable, we have that 
 $$C_{F_n\rtimes_\phi \Z}(t^ax)=\mc S^{-1}\cup   \mc C_{0,a,x,\phi} \cup \mc S$$ is rational and computable.
 
 Since we can test if an element $t^by\in F_n\rtimes_\phi \Z$  belongs to $C_{F_n\rtimes_\phi \Z}(t^ax)$ by checking if $t^byt^ax=t^axt^by$, we can effectively compute a finite set of generators for  $C_{F_n\rtimes_\phi \Z}(t^ax)$.
 \qed\\

\begin{corollary}\label{cor: solseq}
Let $\phi\in \Aut(F_n)$ and $t^au, \, t^bv\in F_n\rtimes_\phi\Z$. The set of solutions  to the equation $x^{-1}(t^au)x=t^bv$ is rational and effectively constructible.
\end{corollary}
\noindent\textit{Proof.} Solving the conjugacy problem in $ F_n\rtimes_\phi\Z$ using \cite[Theorem 1.1]{[BMMV06]}, we check if $t^au$ and $t^bv$ are conjugate. If not, there are no solutions; if they are conjugate, we compute a conjugator $t^cz$ and the set of solutions is $C_{F_n\rtimes_\phi \Z}(t^au)t^cz$, which, by Theorem \ref{thm centralizers}, is computable.
\qed\\

Since the union of two rational subsets is again rational and effectively constructible, it follows from the previous corollary that the set of solutions of  $x^{-1}(t^au)x\in K$ is rational and computable for a finite subset $K$.
Example 5.3 in \cite{[LS11]} shows that the finite subset $K$ cannot be replaced by  arbitrary rational subsets in this statement, as the authors present a [f.g. free]-by-cyclic group $G$, an element $g\in G$ and a rational subset $K\subseteq G$ such that the solution set of $x^{-1}gx\in K$ is not rational.

In a group $G=\langle A\rangle$ we denote by $\pi$ the canonical surjective homomorphism $\pi: (A\cup A^{-1})^*\to G$. Since context-free languages are closed under intersection with rational languages and it is decidable whether a context-free grammar generates the empty language, the following corollary is immediate.
\begin{corollary}\label{cor: cf constraint}
The conjugacy problem with context-free constraints is decidable for [f.g. free]-by-cyclic groups, i.e., there is an algorithm taking as input an two elements $t^ax,\, t^by \in F_n\rtimes_{\phi}\Z$ and a context-free grammar generating a language $L$ such that $L=K\pi^{-1}$ for some $K\subseteq F_n\rtimes_\phi \Z$, that decides whether there is an element $t^cz\in K$ such that $t^ax=(t^cz)^{-1}(t^by) t^c z$. 
\end{corollary}

Naturally the same result holds for subsets $K$ for which we can decide if $K\pi^{-1} \cap L =\emptyset$ for rational languages $L$.

\section{Further work}
We now present the main questions arising from this work.

Firstly, we would also like to understand what $CF(F_n\rtimes\Z)$ consists of. As remarked in the Section 2, for free groups $\Rat(F_n)\subseteq \CF(F_n)$ and for virtually abelian groups $\CF(F_n)\subseteq \Rat(F_n)$. Also for free groups, context-free subsets are well-described. Understanding context-free subsets of [f.g. free]-by-cyclic groups would allow us to fully understand the strength of Corollary \ref{cor: cf constraint}.

\begin{problem}
Can we understand what the class of context-free subsets of [f.g. free]-by-cyclic groups consists of? Is there any relation with the class of rational subsets?
\end{problem}

We also remark that Theorem \ref{thm centralizers} holds for any $G$-by-$\Z$ group as long as $G$ has finitely generated and computable centralizers, decidable Brinkmann's conjugacy problem and decidable twisted conjugacy problem with rational and computable twisted conjugacy classes. Notice that decidability of Brinkmann's conjugacy problem and of the twisted conjugacy problem in $G$ implies decidability of the conjugacy problem in $G$-by-$\Z$ groups, and so the corollaries of the theorem also follow.

If $G$ is a braid group, then the centralizers of its elements are computable \cite{[FG03]} and both Brinkmann's conjugacy problem and the twisted conjugacy problem are decidable \cite{[GV14]}. So, a natural problem is the following:

\begin{problem}
Are twisted conjugacy classes in Braid groups rational and computable?
\end{problem}
We remark that, it follows from the proof of \cite[Theorem 4.9]{[GV14]} that it is enough to solve the problem above for $\varepsilon$-twisted conjugacy classes.

If $G$ is a finitely generated virtually free group, then its centralizers are computable (see for example \cite{[LS11]}) and the twisted conjugacy problem is proved to be decidable in \cite[Theorem 5.4]{[Car23b]}. In the proof, it is shown that, given a virtually free group $G$, an endomorphism $\phi$, elements $u,v\in G$ and constructing $G_1=G*\langle x,y|\rangle$, the free product of $G$ with a free group of rank 2, it is possible to define an endomorphism $\psi$ of $G_1$ such that the set of $\phi$-twisted conjugators of $u$ and $v$  is precisely $x^{-1}\text{Fix}(\psi)y^{-1}\cap G$. The latter is the intersection of two rational and computable (see \cite{[Car22e]} for the computability of Fix($\psi$)) subsets of the virtually free group $G_1$, thus computable \cite[Lemma 4.4]{[Sil02b]}.
The obstruction to a generalization of our Theorem to [f.g. virtually free]-by-$\Z$ is the solution of Brinkmann's conjugacy problem for virtually free groups. We remark that Brinkmann's equality problem was solved in \cite[Theorem 5.3]{[Car23b]}.
\begin{problem}
Is Brinkmann's conjugacy problem decidable for automorphisms of virtually free groups?
\end{problem}

An easier, but still interesting decision problem is satisfiability of the condition in Lemma \ref{centralizerfree}, which can be seen as periodicity modulo conjugation, for some natural classes of groups.
\begin{problem}
Given a f.g. virtually free (resp. hyperbolic) group $G$, an automorphism $\phi\in \Aut(G)$ and an element $x\in G$, can we decide whether there is some $k\neq 0$ such that $x\phi^k\sim x$?
\end{problem}

Another possibility is to extend these results to ascending HNN-extensions of free groups. The proof of the conjugacy problem for [f.g. free]-by-cyclic groups in \cite{[BMMV06]} was generalized to ascending HNN-extensions of free groups in \cite{[Log22]}. An adaptation of our Theorem to this class of groups should be significantly more involved, but might be possible. Also, we remark that the constrained conjugacy problem for ascending HNN-extensions of free groups with finitely gneerated subgroups as constraints is particularly relevant, since \emph{many} 1-relator groups appear as subgroups of ascending HNN-extensions of free groups.

\begin{problem}
Can Theorem \ref{thm centralizers} be extended to ascending HNN-extensions of finitely generated free groups?
\end{problem}

\section*{Acknowledgements}
This work is funded by national funds through the FCT - Fundação para a Ciência e a Tecnologia, I.P., under the scope of the projects UIDB/00297/2020 and UIDP/00297/2020 (Center for Mathematics and Applications).

\bibliographystyle{plain}
\bibliography{Bibliografia}

 \end{document}